\documentstyle[12pt]{article}
\catcode`\@=11
\@addtoreset{equation}{section}

\catcode`\@=12
\newtheorem{Theorem}{Theorem}[section]

\newtheorem{Corollary}{Corollary}[section]
\newtheorem{Note}{Note}[section]
\title{Complex manifolds with generating tangent bundles}

\author{Renyi Ma\\
Department of Mathematical Sciences \\
Tsinghua University \\
Beijing, 100084\\
People's Republic of China\\
rma@math.tsinghua.edu.cn}

\date { }

\begin{document}
\textwidth=125mm
\textheight=185mm
\parindent=8mm
\frenchspacing
\maketitle

\begin{abstract}
Let $M$ be a close complex manifold and $TM$ its holomorphic tangent
bundle. We prove that if the global holomorphic sections of tangent
bundle generate each fibre, then $M$ is a complex homogeneous
manifold. Our proof depends on the complex version of
Chow-Rashevskii theorem in Carnot-Caratheodory spaces.
\end{abstract}

\section{Introduction and results }

The main result is the following theorem.

\begin{Theorem}
Let $M$ be a close complex manifold and $TM$ its holomorphic
tangent bundle. We prove that if the global holomorphic sections
of tangent bundle generate each fibre, then $M$ is isomorphic to a
complex homogeneous manifold.
\end{Theorem}

It is obvious that Theorem1.1 implies

\begin{Corollary}
Every compact K\"ahler manifold with Griffith's ample tangent bundle
(see\cite{gri}) is isomorphic to the projective space.
\end{Corollary}

 {\bf Sketch of proof on
Theorem 1.1.} In section 2, we give a complex version of
Chow-Rashevskii theorem in Carnot-Caratheodory geometry. In Section
3, we give the proof of Theorem1.1.

\section{Complex Chow-Rashevskii Theorem}

In this section, we prove the complex version of Chow-Rashevskii's
connectivity theorem.
\begin{Theorem}
Let $M$ be a close connected complex manifold and $TM$ its
holomorphic tangent bundle. If the global holomorphic sections of
tangent bundle generate each fibre, then for any two points
$p,q\in M$, there exists a holomorphic automorphism $T:M\to M$
such that $T(p)=q$ and $T$ is isotopic to $Id:M\to M$ in the group
of automorphism of $M$.
\end{Theorem}
Proof. Since $M$ is a complex close manifold, the complex vector
space $V$ generated by the holomorphic vector fields has finite
dimension $p$. We assume that the holomorphic vector fields
$Y_1,....,Y_p$ generate $V$. Since $M$ is close, the vector field
$Y_i$ integrates a holomorphic flow $Y_i(z),z\in C$, $i=1,...,p$.
For each $m\in M$, we consider the composed action map $E_m:C^p\to
M$ defined by
\begin{equation}
(z_1,...,z_n)\to Y_1(z_1)\circ ...\circ Y_p(z_p)(m).
\end{equation}
The differential of $E_m$ at the origin $0\in C^p$ sends the $C^p$
onto the span of the fields $Y_i$ in $T_m(M)$ and, hence, is
surjective by assumption of the theorem. Thus the orbit $G(m)$ is
open in $M$ for each $m\in M$ by the implicit function
theorem(see\cite{gro}). It is easy to see that $G(m)$ is close in
$M$. Since $M$ is connected, $G(m)=M$, q.e.d.

\begin{Note}
Our proof is inspired by Gromov's proof on Chow-Rashevskii theorem
in Carnot-Caratheodory geometry(see\cite{gro}).
\end{Note}

\section{Proof on Theorem 1.1}

{\bf Proof of Theorem1.1:} by Theorem2.1, it is obvious that the
complex Lie group $G$, i.e., the automorphism group of $M$ acts
transitively on $M$. So, $M$ is complex homogeneous manifolds. This
yields Theorem1.1.

\end{document}